\documentclass[twocolumn,twoside]{IEEEtran}
\usepackage{times}
\usepackage{amssymb}
\usepackage{theorem}
\usepackage{amstext}    
\usepackage{mathrsfs}   
\usepackage{psfig}

\renewcommand{\leftmark}%
{{\sc IEEE Transactions on Information Theory}, to appear}
\renewcommand{\rightmark}%
{Submitted August 12, 2003; revised March 28, 2004}

\renewcommand{\Bbb}{\mathbb}

\newcommand{\Fn}{{\Bbb{F}}_{\hspace{-1pt}2}^{\hspace{1pt}n}}
\newcommand{\Fnw}{{\Bbb{F}}_{\hspace{-1pt}2}(n,w)}
\newcommand{\R}{{\Bbb{R}}}
\newcommand{\C}{{\Bbb{C}}}

\newcommand{\GG}{\cG_{\rm G}}
\newcommand{\GS}{\cG_{\cS}}
\newcommand{\GGq}{\cG_{q,{\rm G}}}
\newcommand{\GSq}{\cG_{q,{\cS}}}

\newcommand{\be}[1]{\begin{equation}\label{#1}}
\newcommand{\ee}{\end{equation}} 
\newcommand{\eq}[1]{(\ref{#1})}

\newcommand{\Proof}{\hspace*{5pt}{\em Proof.}\hspace*{1.0ex}}
\def\qed{\hskip 3pt \hbox{\vrule width4pt depth2pt height6pt}}

\newcommand{\cA}{{\cal A}}
\newcommand{\cB}{{\cal B}}

\newcommand{\cE}{{\cal E}}

\newcommand{\cG}{{\cal G}}
\newcommand{\cH}{{\cal H}}
\newcommand{\cI}{{\cal I}}

\newcommand{\cS}{{\cal S}}

\newcommand{\Tref}[1]{Theorem\,\ref{#1}}
\newcommand{\Pref}[1]{Proposition\,\ref{#1}}
\newcommand{\Lref}[1]{Lemma\,\ref{#1}}

\renewcommand{\le}{\leqslant}
\renewcommand{\leq}{\leqslant}
\renewcommand{\ge}{\geqslant}
\renewcommand{\geq}{\geqslant}

\newcommand{\dfn}{\it}

\newcommand{\deff}{\mbox{$\stackrel{\rm def}{=}$}}

\newcommand{\Strut}[2]{\rule[-#2]{0cm}{#1}}

\newcommand{\trunc}[1]{\left\lfloor {#1} \right\rfloor}
\newcommand{\ceil}[1]{\left\lceil {#1} \right\rceil}

\newcommand{\sfrac}[2]{\mbox{\raisebox{.8mm}{\footnotesize $\scriptstyle #1$} 
                          \footnotesize$\!\!\! / \!\!$
                          \raisebox{-.8mm}{\footnotesize $\scriptstyle #2$}}}

\newcommand{\wt}{{\rm wt}}

\newcommand{\zero}{{\mathbf 0}}
\newcommand{\one}{{\mathbf 1}}

\newcommand{\al}{\alpha}

\newcommand{\GV}{Gil\-bert-Var\-sha\-mov}
\newcommand{\And}{A_2(n,d)}
\newcommand{\Andw}{A(n,2d,w)}

\newcommand{\fGV}{f_{\rm GV}(n,d)}
\newcommand{\fV}{f_{\rm V}(n,d)}
\newcommand{\fE}{f_{\rm E}(n,d)}
\newcommand{\fT}{f_{\rm T}(n,d)}
\newcommand{\fBGS}{f_{\rm BGS}(n,d)}
\newcommand{\fFone}{f_{{\rm F}_1}(n,d)}
\newcommand{\fFtwo}{f_{{\rm F}_2}(n,d)}

\theoremstyle{plain} 
\theorembodyfont{\normalfont\itshape}

\newtheorem{thm}{Theorem} 
\newenvironment{theorem}
{\begin{thm}\hspace*{-1ex}{\bf.}}{\end{thm}}

\newtheorem{prop}[thm]{Proposition} 
\newenvironment{proposition}
{\begin{prop}\hspace*{-1ex}{\bf.}}{\end{prop}}

\newtheorem{lem}[thm]{Lemma}
\newenvironment{lemma}{\begin{lem}\hspace*{-1.0ex}{\bf.}}{\end{lem}}

\newtheorem{cor}[thm]{Corollary}

\newtheorem{defi}{Definition}

\theorembodyfont{\normalfont}

\newtheorem{exam}{Example}

\outer\def\proclaim #1. #2\par{\medbreak
 \noindent{\bf#1.\enspace}{\sl#2\par}%
 \ifdim\lastskip<\medskipamount \removelastskip\penalty55\medskip\fi}

\mathchardef\inn="3232
\renewcommand{\in}{\mbox{$\,\inn\,$}}

\makeatletter
\renewcommand{\@endtheorem}{\endtrivlist}
\makeatother

\newcounter{enumrom}
\renewcommand{\theenumrom}{(\roman{enumrom})}

\renewcommand{\thebibliography}[1]{\section*{References}\footnotesize\list
    {[\arabic{enumi}]\hfill}{\settowidth\labelwidth{[#1]}\leftmargin\labelwidth
    \advance\leftmargin\labelsep \itemsep 0pt plus .5pt
    \usecounter{enumi}}
    \def\newblock{\hskip .11em plus .33em minus .07em}
    \sloppy\clubpenalty4000\widowpenalty4000
    \sfcode`\.=1000\relax}

\makeatletter

\gdef\@punct{.\ \ }  
\def\@sect#1#2#3#4#5#6[#7]#8{%
  \ifnum #2>\c@secnumdepth
     \def\@svsec{}
  \else
     \refstepcounter{#1}\edef\@svsec{%
     \ifnum #2>0{{\csname the#1\endcsname}}.\fi%
    \hskip .5em}
  \fi
  \@tempskipa #5\relax
  \ifdim \@tempskipa>\z@
     \begingroup #6\relax
       \@hangfrom{\hskip #3\relax\@svsec}{\interlinepenalty \@M #8\par}
     \endgroup
     \csname #1mark\endcsname{#7}
     \addcontentsline{toc}{#1}{\ifnum #2>\c@secnumdepth\else
          \protect\numberline{\csname the#1\endcsname}\fi#7}
  \else
     \def\@svsechd{#6\hskip #3\@svsec #8\@punct\csname #1mark\endcsname{#7}
     \addcontentsline{toc}{#1}{\ifnum #2>\c@secnumdepth \else
          \protect\numberline{\csname the#1\endcsname}\fi#7}}
  \fi
  \@xsect{#5}}

\def\@ssect#1#2#3#4#5{\@tempskipa #3\relax
  \ifdim \@tempskipa>\z@
    \begingroup #4\@hangfrom{\hskip #1}{\interlinepenalty \@M #5\par}\endgroup
  \else \def\@svsechd{#4\hskip #1\relax #5\@punct}\fi
  \@xsect{#3}}

\makeatother

\title{\mbox{Asymptotic Improvement of the~Gilbert-Varshamov}\\ Bound
on the Size of Binary Codes\\[.50ex]}

\author{{\Large\sc Tao Jiang} 
        {\large ~and~} 
        {\Large\sc Alexander Vardy}\vspace{2ex}
\thanks{Manuscript submitted August\,12,\,2003, and \looseness=-1
revised \mbox{March\,28,\,2004}.
This work was supported in part 
by the David and Lucile Packard Foundation, by 
the National Science Foundation, and by the Miami 
University Summer Faculty Research Grant.
The material~in~this~paper was presented in part at 
the 41-st Annual Allerton Conference on 
Communication, Control, and Computing,
Mon\-ticello, IL, October 2003.
}
\thanks{Tao Jiang is with the
Department of Mathematics and Statistics,
Miami University, 
Oxford, OH 45056, U.S.A.
(e-mail: {\tt jiangt@muohio.edu}).}
\thanks{Alexander Vardy is with 
the Department of Electrical and Compu\-ter Engineering,
the Department of Computer Science and Engineering, 
and the Department of Mathematics,
University of California San Diego,
La Jolla, CA 92093-0407, U.S.A.
(e-mail: {\tt vardy@kilimanjaro.ucsd.edu}.)}
\thanks{Communicated by Ralf Koetter, Associate Editor for Coding Theory.}
\vspace*{-.25cm}}

\begin{document}

\maketitle

\begin{abstract}
\looseness=-1
Given positive integers $n$ and $d$, 
let $A_2(n,d)$ denote the maximum size
of a binary code of length $n$ and minimum distance $d$. 
The well-known Gilbert-Varshamov bound asserts that 
$A_2(n,d) \geq 2^n/V(n,d{-}1)$,
where $V(n,d)=\sum_{i=0}^{d}\!{n\choose i}$
is the volume of a Hamming sphere of radius~$d$.
We show that, in fact, 
there exists a~positive constant $c$ such that
$$
A_2(n,d) \ \geq\  c\: \frac{2^n}{V(n,d{-}1)} \, \log_2\! V(n,d{-}1)
$$
whenever $d/n \le 0.499$.
The result follows by recasting the Gilbert-Var\-shamov bound
into a graph-theoretic framework and using the fact that the
corresponding graph is locally sparse. Generalizations and 
extensions of this result are briefly discussed.\vspace{1ex}
\end{abstract}


\section{Introduction}  
\label{sec1}

\noindent \looseness=-1
Let $A_q(n,d)$ denote the maximum number of codewords in
a~code of length $n$ and
minimum distance $d$ over an alphabet with $q$ letters. The \GV\ 
bound, which asserts that
\be{q-GV}
A_q(n,d)
\ \ge \
\frac{q^n}{\sum_{i=0}^{d-1} {n \choose i} (q\,{-}\,1)^i}
\ee
is one of the most well-known and fundamental results in coding 
theory. In this paper, we focus on binary codes (although an extension
of our results to codes over an arbitrary alphabet is 
discussed in Section\,V). Thus we let
$$
V(n,d) 
\ \ \deff \ \
\sum_{i=0}^{d} {n \choose i}
$$
de\-note the volume of a Hamming sphere of radius $d$ in $\Fn$,
and~consider the binary version of~\eq{q-GV}, namely
\be{2-GV}
A_2(n,d)
\ \ge \
f_{\rm GV}(n,d) \ \ \deff\ \ \frac{2^n}{V(n,d{-}1)}
\ee
This inequality was first proved by Gilbert~\cite{Gilbert} in 1952.
It was subsequently improved upon by Varshamov~\cite{Varshamov}.
However, following the established terminology, we will
refer to \eq{q-GV}~and~\eq{2-GV} as the Gilbert-Varshamov 
bound.
This bound is used extensively in the coding theory 
literature~\cite{MWS,Handbook}, and has been generalized
to numerous contexts~\cite{CLH,Levenshtein71,GF,MR,Tolhuizen}.

\looseness=-1
Improving upon the \GV\ bound asymptotically is 
a notoriously difficult task~\cite{Handbook,TV}.
The breakthrough work of 
Tsfasman-Vl\v{a}du\c{t}-Zink~\cite{TVZ}
led to an asymptotic improvement of~\eq{q-GV},~but only for
alphabets of size $q \ge 49$ 
(see also the recent papers~\cite{Elkies,Xing}). 
For $q < 46$, no asymptotic
improvements upon~\eq{q-GV} 
are currently known~\cite{ZS}.
In fact, a well-known conjecture (cf.~Goppa~\cite{Goppa})
asserts that the binary version~\eq{2-GV} of the \GV\
bound is asymptotically exact.

\looseness=-1
Nevertheless, for small $n$ and $d$, the size of best known
binary codes~\cite[Chapter\,5]{Handbook}
often exceeds $\fGV$ 
by a large factor. Thus it is natural to ask whether the bound~\eq{2-GV} 
can be strengthened. Indeed, various improvements
upon the binary \GV\ bound were presented (in chronological order)
by 
Varshamov~\cite{Varshamov}, 
Hashim \cite{Hashim}, 
Elia~\cite{Elia},
Tolhuizen~\cite{Tolhuizen}, 
Barg-Guritman-Simonis~\cite{BGS}, and
Fabris~\cite{Fabris}. 
We review these improvements in detail
in the next section.
One of our main results herein is the following
theorem, which strengthens the \GV\ bound using a technique quite
different from those 
of~\cite{BGS,Elia,Fabris,Hashim,Tolhuizen},
and \cite{Varshamov}.\vspace{-0.50ex}
\begin{theorem}
\label{exact}
For $x \in \R$, let\/ $\ceil{x}^+$ 
denote the smallest {nonnegative} integer $m$ with~$m \ge x$.
Given positive integers $n$ and $d$, with $d \le n$, let 
$e(n,d)$ denote the following quantity
$$
e(n,d) 
\ \ \deff\ \
\frac{1}{6}
\sum_{w=1}^{d} \hspace{-.5ex} {n \choose w}\hspace{-.5ex}
\left(
\sum_{i=1}^{d}\hspace{-.25ex}
\sum_{j=\ceil{\frac{w+i-d}{2}}^+}^{\min\{w,i\}}\hspace{-1ex}
{w \choose j}\! {n{-}w \choose i\,{-}\,j} - 1\hspace{-.5ex}
\right) 
$$
Then 
\be{th1}
A_2(n,d)
\ \ge \
\frac{2^n}{V(n,d{-}1)} 
\cdot
\frac{\log_2\! V(n,d{-}1)\ - \ \log_2 \sqrt{e(n,d{-}1)}}{10}\vspace{1ex}
\ee
\end{theorem}

What distinguishes \Tref{exact} from prior improvements of the \GV\
bound is the asymptotic behavior of~\eq{th1}. 
All the previously known 
explicit lower bounds on $\And$ that we are aware of, including those
of~\cite{BGS,Elia,Fabris,Hashim,Tolhuizen}, and~\cite{Varshamov}, 
have the following
property: if we write the bound as \mbox{$\And \ge f(n,d)$}, then
\be{f(n,d)}
f(n,d) 
\ = \ 
O\Bigl(\fGV\Bigr)
\ee
In fact, as we shall see in the next section,
for some of these bounds 
$
f(n,d) = \fGV \left(\Strut{2ex}{0.1ex} 1 + o(1) \right)
$, 
where $o(1)$ tends to zero exponentially fast with $n$.
In contrast, the asymptotic behavior of~\eq{th1} 
is given by the following theorem.\pagebreak[3.99]

\begin{theorem}
\label{asymptotic}
Let $n$ and $d$ be positive integers, with $d/n \le 0.499$. Then 
there exists~a~positive constant $c$ such that
\be{th2}
A_2(n,d) \ \geq\  c\: \frac{2^n}{V(n,d{-}1)} \, \log_2\! V(n,d{-}1)
\vspace{1.25ex}
\ee
\end{theorem}

\noindent
{\bf Remark.}
The constant in \Tref{asymptotic}, the way it is 
stated~above, may depend on the ratio $d/n$. However,
if we only wish to claim that \eq{th2} is true for
all {sufficiently large $n$}, then $c$ becomes
an absolute constant, independent of both $n$ and $d$.
For more on this, see~\eq{th2-proof}.
Also note that while the bound in~\eq{th2}
holds for any $n$ and $d$ with $d/n \le 0.499$,
it is useful only when the ratio $d/n$ is constant. 
If we allow $d/n \to 0$ as $n \to \infty$, then better 
bounds on $\And$ are known~\cite{BJ,MWS,TV}.
\vspace{1ex}

\looseness=-1
So, how does \Tref{asymptotic} relate to the conjecture that
the \GV~bound is asymptotically exact for $q=2$?
This depends on the interpretation. 
If one views the conjecture as dealing with the asymptotics
of $\And$ itself, namely 
the \underline{size} of the best binary codes, then it
corresponds to the assertion that for all positive 
$\delta < 0.5$, we have
\be{GV-conjecture}
\lim_{n \to \infty} {A_2(n,\delta n) \over f_{\rm GV}(n,\delta n)}
\ = \
\text{const}
\ee
where the constant might be a function of $\delta$.
\mbox{\Tref{asymptotic} clearly}
shows that this is false:~the limit
$\lim_{n \to \infty} A_2(n,\delta n)/f_{\rm GV}(n,\delta n)$
does not exist for any $\delta$.
Indeed, the theorem implies that
\be{GV-revised}
\log_2 \And
\ \ge \
\log_2 \fGV \ + \ \log(n) \ + \ \text{const} \ + \ o(1)
\ee
On the other hand, 
it is more common to interpret the conjecture as dealing
with the asymptotics of the best possible \underline{rate} of 
a binary code, namely the function $R(n,d) = \log_2 \And/n$.
In this case, the conjecture
could still be true, since 
the term $\log(n)/n$ will vanish for $n \to \infty$.

The rest of this paper is organized as follows. In the next section,
we review the previously known improvements of the \GV\ bound,
with the aim of establishing~\eq{f(n,d)}. In~Section\,III, we recast
the problem of estimating $\And$ into a~graph-theoretic framework,
and express $\And$ as the independence number of a certain graph
(\Lref{lemma-Gilbert}). We then recover the \GV\ bound as 
a straightforward consequence of a simple bound on the
independence number of a~graph (\Pref{P-GV}).
The key idea in the proof of Theorems \ref{exact} and \ref{asymptotic}
is surprisingly simple: the bound on the independence number 
used in \Pref{P-GV} 
can be improved upon, providing 
the graph at hand is locally sparse (\Tref{regular}). 
In Section\,IV, we show that this is, indeed, the case.
Specifically, we derive a simple closed-form expression 
for the number of edges in the relevant graph (\Pref{exact-e(GS)}),
and then prove that this graph is sparse for all sufficiently
large $n$ whenever $d/n \le 0.4994$ (\Pref{P-asymptotic}). This 
completes the proof of Theorems~\ref{exact}~and~\ref{asymptotic}.
In Section\,V, we briefly discuss various extensions and generalizations
of our results. In particular, we show that just like the bounds 
of Gilbert~\cite{Gilbert} and Varshamov~\cite{Varshamov}, our bound
can be proved ``constructively.'' That is, there is an (exponential-time)
algorithm~\cite{HL} that actually constructs codes 
satisfying~\eq{th2}. We also generalize \Tref{exact} to arbitrary
alphabets (\Tref{q-ary}) and to constant-weight codes. 
Finally, we point out a number of intriguing
open problems related to the results of this paper.\pagebreak[3.99]

\vspace{3ex}
\section{Comparison with Prior Work}  
\label{sec2}

\noindent\looseness=-1
In this section, we briefly review previously known (to us) improvements 
of the \GV\ bound~\eq{2-GV}, roughly in chronological order, and 
establish the claim of~\eq{f(n,d)}.

The first improvement on~\eq{2-GV} is due to Varshamov himself. 
Varshamov showed in~\cite{Varshamov}~that $\And \ge \fV$,
where\footnote
{Usually, 
$\fV$ is defined as $2^k$, where $k$ is the largest integer satisfying 
$2^k < 2^n/V(n{-}1,d{-}2)$. The
explicit form~\eq{Varshamov-f} 
is equivalent to this definition.\vspace*{-.75ex}
}\be{Varshamov-f}
f_{\rm V}(n,d) \ \ \deff\ \ 
\frac{2^{n-1}}{\displaystyle 2^{\trunc{\log_2\! V(n-1,d-2)}}}
\ee
and, moreover, there exist {linear\/} codes that attain this
bound. We now show that the ratio $\fV/\fGV$ is upper bounded 
by a constant. Indeed, we have
$$
\frac{\fV}{\fGV}
\ \le \
\frac{V(n,d{-}1)}{V(n{-}1,d{-}2)}
\ = \ 
1 \ + \ \frac{V(n{-}1,d{-}1)}{V(n{-}1,d{-}2)}
$$
where the equality above follows 
from~the~fact~that~$V(n,d{-}1) = V(n{-}1,d{-}1) +  V(n{-}1,d{-}2)$. 
Expressing $V(n,d)$ as the sum $\sum_{i=0}^{d}\!{n\choose i}$,
we further obtain\vspace{-.25ex}
\begin{eqnarray} %
\frac{V(n{-}1,d{-}1)}{V(n{-}1,d{-}2)}
& \!{=}\! &
1 +\, \frac{{n{-}1 \choose d{-}1}}{V(n{-}1,d{-}2)}\\
\label{Varshamov-auxA}
& \!{=}\! &
1 +\, \frac{1}{\sum_{i=0}^{d-2} \frac{ (d-1)!\,\,(n-d)! }{ i!\ (n-i-1)!}}\\
\label{Varshamov-aux}
& \!{\le}\! &
1 +\, \frac{n - (d{-}1)}{d-1}
\end{eqnarray}
where the inequality in \eq{Varshamov-aux}
follows by retaining a single term in the sum of \eq{Varshamov-auxA}, 
namely the term  
corresponding to $i = d{-}2$.
Thus $\fV/\fGV \le (\delta + 1)/\delta$, where $\delta = (d{-}1)/n$.

Another improvement of~\eq{2-GV} was proposed 
by Hashim~in~1978.  
Hashim~\cite[eq.\,(7)]{Hashim} proved the following. 
Let $t = \ceil{(d{-}1)/2}$ and let $A(w;n,k,d)$ denote the minimum 
number of codewords of weight $w$ in an $(n,k,d)$ binary linear code.
Then~\cite{Hashim} shows that $\And \ge 2^k$, where $k$ is the largest 
integer satisfying\vspace{-.25ex}
\be{Hashim}
V(n{-}1,d') 
\ - \sum_{w=d}^{d'+\,t} 
 \sum_{i=w-d'}^t \hspace{-.75ex} {w \choose i} A(w;n,k,d)
\ < \ 
2^{n-k}
\ee
where $d' = d-2$. Unfortunately, the bound \eq{Hashim}
is non-explicit. Hashim~\cite{Hashim} writes
that ``this improved bound requires the determination of the lowest
possible value of $A(w;n,k,d)$, where $w = d,d{+}1,\ldots,d{-}2{+}t$,
in terms of the code parameters $n$, $k$, and $d$.''
While various estimates of $A(w;n,k,d)$ are known \cite{ABL,KL,KL1,KL2},
we are not aware of any results that can be used in conjunction 
with~\eq{Hashim} to produce an explicit lower bound on $\And$,
at least not without a substantial research effort.

\looseness=-1
In 1983, 
Elia\,\cite{Elia} has extended the Varshamov bound~\eq{Varshamov-f}
in a~different~way. Specifically, it is shown~in 
\cite[Corollary\,2]{Elia} that
$\And \ge \fE$, where\vspace{-.50ex}
$$
\fE \ \ \deff\ \ 
\frac{2^{n-2}}
{\displaystyle 
\max\left\{2^{\trunc{\log_2\! V(n-3,d-2)}},\ %
           2^{\trunc{\log_2\! V(n-2,d-3)}}%
\right\}}
$$
It is not difficult to see that, again, 
the ratio $\fE/\fGV$ is upper bounded\pagebreak[3.99]
by a~constant. Indeed, writing $2^{\trunc{\log_2\! V(n,d)}}$
as $V(n,d)/2^{\{\log_2\! V(n,d)\}}$, where 
$\{x\}$ denotes the fractional
part of $x\,{\in}\,\R$, we have
$$
\fE 
\ \le \
\frac{2^{n-2}2^{\{\log_2\! V(n-3,d-2)\}}}{V(n{-}3,d{-}2)}
\ \le \
\frac{2^{n-1}}{V(n{-}3,d{-}2)}
$$
This, in turn, leads to the following bound
\be{Elia-aux}
\frac{\fE}{\fGV}
\ \le \
\frac{V(n,d{-}1)}{2V(n{-}3,d{-}2)}
\ \le \
\frac{8 V(n{-}3,d{-}1)}{2 V(n{-}3,d{-}2)}
\ee
We know from~\eq{Varshamov-aux} that 
$
{V(n{-}3,d{-}1)}/{V(n{-}3,d{-}2)} \le {1}/{\delta}
$,
where $\delta = (d{-}1)/n$.
In conjunction with~\eq{Elia-aux}, this implies
that $\fE/\fGV \le 4/\delta$. 

Tolhuizen~\cite{Tolhuizen} established yet another slight improvement
of \eq{2-GV} using Tur\'{a}n's theorem~\cite[Chapter\,4]{VW}.
Specifically, Tolhuizen~\cite[p.\,1605]{Tolhuizen} shows that
$\And \ge \fT + 1$, where $\fT$ is the largest integer satisfying
\be{Tolhuizen}
\frac{2^n}{\fT} \ + \ \frac{r(\fT - r)}{2^n \fT}
\ > \
V(n,d{-}1)
\ee
with $r$ being the remainder when $2^n$ is divided by $\fT$.
If we ignore the second term on the left-hand side of~\eq{Tolhuizen},
then this is precisely the \GV\ bound~\eq{2-GV}. Otherwise, it is easy
to see that\vspace{-1ex}
\begin{eqnarray*}
\fT
& \hspace{-1ex}\le\hspace{-1ex} &
\frac{2^n}{V(n,d{-}1) - 2^{-(n+2)}}\\
& \hspace{-1ex}\le\hspace{-1ex} &
\frac{2^n}{V(n,d{-}1)}
\cdot 
\frac{2^{n+2}}{2^{n+2} - 1}
\ = \
\fGV \left(\Strut{2ex}{0.1ex} 1 + o(1) \right)
\end{eqnarray*}

The latest improvement on \eq{2-GV} 
is due to Fabris~\cite{Fabris}. In fact, Fabris~\cite{Fabris}
proves~two new bounds on $\And$. The first bound is given by
$\And \ge \fFone$, where
\be{Fabris1}
\fFone \ \ \deff\ \ 
\frac{2^n - \cI(n,d{-}1)}{V(n,d{-}1) - \cI(n,d{-}1)}
\ee
and $\cI(n,d{-}w)$ is the volume of the intersection between 
two Hamming spheres of radius $d-w$, whose centers are  
distance $d$ apart. 
The second bound is $\And \ge \fFtwo$,~where
\be{Fabris2}
\fFtwo \ \ \deff\ \ 
\frac{2^n}{V(n,d{-}1)}
\left(
\frac{V(n,d{-}1) +\cI(n,d{-}2)}{V(n,d{-}2)}
\right)
\ee
Obviously $\cI(n,d{-}2) \le V(n,d{-}2)$. Thus it follows 
straightforwardly from~\eq{Fabris2},\,\eq{Varshamov-aux}
that 
$\fFtwo/\fGV \le (\delta + 1)/\delta$. It is not 
difficult to see (cf.~\Lref{exact-degree}) that
$$
\cI(n,d-w)
\ = \
\sum_{i=w}^{d-w}
\sum_{j=\ceil{\frac{w+i}{2}}}^{i}
{w \choose j} {n\,{-}\,w \choose i\,{-}\,j}
$$
In Section\,IV herein, we will show (in a different context) that 
\mbox{$\lim_{n \to \infty} \cI(n,d{-}1)/V(n,d{-}1) =0$}. 
In conjunction with \eq{Fabris1},\linebreak[3.99]
this immediately implies that
\mbox{$\fFone = \fGV \! \left(\Strut{2ex}{0.1ex} 1 + o(1) \right)$}.

Finally, 
the recent work of Barg, Guritman and Simonis~\cite{BGS} contains various
extensions and generalizations of the Varshamov bound~\eq{Varshamov-f}
as well as related prior work by Hashim~\cite{Hashim}, Elia~\cite{Elia},
and Edel~\cite{Edel}. However, just as the Hashim~bound, 
most of the results~of~\cite{BGS} are non-explicit --- they provide
methods for constructing codes, 
but a substantial research effort would be required to convert them into 
an explicit lower bound on $\And$.
On the other hand, \cite{BGS} does contain the
following generalization of Elia's bound: for all 
$b = 0,1,\ldots,d{-}1$, if
\mbox{$2^{b-1}V(n{-}b,d{-}b{-}1) < 2^{n-k}\!$ and there exists 
an $(n{-}b,\kern-1pt k{-}1,\kern-1pt d)$} 
code, then $\And \ge 2^k$.
If we use the Varshamov bound~\eq{Varshamov-f} to guarantee the
existence of the $(n{-}b,k{-}1,d)$ code, then this reduces to
$\And \ge \fBGS$, where\vspace{-.25ex}
\begin{eqnarray*}
\lefteqn{\fBGS \ \ \ \deff\ \ }\\[-.5ex]
& & \frac{2^{n}}
{\displaystyle 
2^b \max\left\{2^{\trunc{\log_2\! V(n-b-1,d-2)}},\ %
           2^{\trunc{\log_2\! V(n-b,d-b-1)}}%
\right\}}
\end{eqnarray*}
with $b$ serving as an optimization parameter (note that for $b=1$,
we recover the Varsha\-mov bound~\eq{Varshamov-f}, while for $b=2$
this is precisely the Elia bound). 
Proceeding as in~\eq{Varshamov-aux} and \eq{Elia-aux}
while taking into account that 
$V(n,d{-}1) \le 2^{b+1} V(n{-}b{-}1,d{-}1)$,~it~is 
easy to see that $\fBGS/\fGV \le 4/\delta$.

\vspace{2.5ex}
\section{Gilbert-Varshamov Bound\\ and Locally Sparse Graphs}
\label{sec3}

\noindent\looseness=-1
We first recall some elementary terminology from graph theory.
A {\dfn graph\/} $G$ consists of a~set of {\dfn vertices} $V(G)$ 
and a set $E(G)$ of pairs of vertices, whose elements are called {\dfn edges}. 
We henceforth assume that both $V(G)$ and $E(G)$ are finite sets.
We use $n(G)$ and $e(G)$ to denote, respectively, 
the number of vertices and
the number of edges in $G$.
Two vertices $u,v \in V(G)$ are {\dfn adjacent} or {\dfn neighbors}
in $G$ iff $\{u,v\} \in E(G)$. The set of all neighbors of 
a~vertex $v$ is denoted $N(v)$ and called the 
{\dfn neighborhood\/} of $v$. The {\dfn degree\/} of a vertex $v \in V(G)$,
denoted $\deg(v)$,~is~defined as $\deg(v) = |N(v)|$. 
A graph $G$ is said to be~\mbox{\dfn $\Delta$-regular\/} if
$\deg(v) = \Delta$ for all $v \in V(G)$.
A set $K \subseteq V(G)$ is a {\dfn clique\/} if every vertex in $K$
is adjacent to all other vertices in $K$. A clique consisting of\, $3$
vertices is a {\dfn triangle}. A~set $\cI \subseteq V(G)$
such~that no two vertices in $\cI$ are adjacent is an 
{\dfn independent set}. 
A {\dfn proper $c$-coloring\/} of $G$ is a~partition of $V(G)$ into $c$ 
independent sets.
The maximum number of vertices 
in an independent set is called the 
{\dfn independence number\/} of $G$, and denoted $\alpha(G)$. 

\looseness=-1
The $n$-dimensional hypercube $\cH_n$ is defined as a graph whose vertex set
$V(\cH_n)$ is the set of all~binary vectors of length $n$, with
$u,v \in V(\cH_n)$ being adjacent iff $d(u,v) = 1$, where $d(\cdot,\cdot)$
is the Hamming distance. 
Note that the graph distance in $\cH_n$ is equal to the Hamming distance.
Given a minimum distance $d$, 
we define the {\dfn Gilbert graph\/} as $\cH_n$ to the power $(d-1)$.

\proclaim Definition.
Let $n$ and $d \le n$ be positive integers.
The corresp\-onding \emph{Gilbert graph $\GG$} is
defined as follows: $V(\GG) = \Fn$ 
and 
$\{u,v\} \in E(\GG)$ if and only if\, $1 \le d(u,v) \le d-1$.

Clearly, a binary code of length $n$ and minimum distance $d$ 
is an independent set in the Gilbert graph $\GG$. 
Conversely, any independent set in $\GG$ is a binary code of length $n$ 
and minimum distance at least $d$. This proves the 
following.\vspace{-1ex}

\begin{lemma}
\label{lemma-Gilbert}\vspace{-1ex}
\be{alpha(GG)}
\And \ = \ \alpha(\GG)
\ee
\end{lemma}

\looseness=-1
\Lref{lemma-Gilbert} makes it possible to recover the \GV\
bound \eq{2-GV} as a straightforward corollary to\pagebreak[3.99] 
a~simple~bound
on the independence number of a graph. Since numerous distinct
proofs of the GV bound (e.g.\ using Tur\'{a}n's 
theorem~\cite{BGS,Tolhuizen}) abound in the literature,
it is somewhat surprising that the simple proof below seems to have
not been previously published.\vspace{-.5ex}

\begin{proposition}
\label{P-GV}\vspace*{-1.25ex}
\be{P-alpha}
\alpha(\GG)
\ \ge \
\frac{2^n}{V(n,d{-}1)}
\ee
\end{proposition}

\Proof 
\looseness=-1
By definition, the Gilbert graph $\GG$ is $\Delta$-regular with 
$\Delta = V(n,d{-}1) - 1$. 
Let $\cI$ be a \emph{maximal} independent set in $\GG$, 
and let $\cE \subset E(\GG)$ be the set of edges
with one endpoint in $\cI$ and the other in $V(\GG) - \cI$.
Since $\cI$ is an independet set, we have $|\cE| = \Delta|\cI|$.
Since $\cI$ is maximal, every vertex of $V(\GG) - \cI$ is
adjacent to at least one vertex of $\cI$, and so 
$|\cE| \ge n(\GG) - |\cI|$. 
Therefore
$
\al(\GG) \ge |\cI| \ge n(\GG)/(\Delta{+}1) = 2^n/V(n,d{-}1)
$.~\qed\vspace{1.5ex}

\noindent
{\bf Remark.}
The trivial bound $\al(\GG) \ge n(\GG)/(\Delta{+}1)$
proved in \Pref{P-GV} is well known in graph theory.
This bound can be strengthened somewhat using 
Brooks' theorem~\cite[p.\,20]{Brooks,VW}:
since $\GG$ is obviously neither a complete 
graph nor an odd cycle, it must be $\Delta$-colorable. 
The largest color class in a proper $\Delta$-co- loring 
of $\GG$ has to contain at least $n(\GG)/\Delta$ vertices. 
\vspace{1ex}

Note that the proof of~\eq{P-alpha} requires very little
information about $\GG$. Thus we can easily improve upon~\eq{P-alpha}
using the fact that the neighborhood $N(v)$ of every vertex $v$ in 
$\GG$ is fairly sparse.
First, we need a couple of well-known results about locally sparse graphs.
We say that $G$ is a graph with maximum
degree at most $\Delta$ if $\deg(v) \le \Delta$ for all $v \in V(G)$. 
\vspace{-1.25ex}

\begin{lemma} 
\label{triangle}
Let $G$ be a graph with maximum degree at most $\Delta$, 
and suppose that there 
are no triangles in $G$. 
Then 
\be{AKS}
\alpha(G)
\ \geq \ 
\frac{n(G)}{8\Delta}\, \log_2 \Delta
\ee
\end{lemma}

\Lref{triangle} was first proved by 
Ajtai, Koml\'os, and Szemer\'e\-di~\cite{AKS}
(but see \cite[p.\,272]{AS} for 
a much shorter proof of the same result). 
Subsequently, the bound in~\eq{AKS} has been
extended from graphs without triangles to 
graphs with relatively few triangles.
In particular, a proof of the following lemma
can be found, for example, 
in Bollob\'as~\cite[Lemma\,15, p.\,296]{Bollobas}.
\vspace{-1.25ex}

\begin{lemma}
\label{sparse}
Let $G$ be a graph with maximum degree at most $\Delta$ 
and suppose that~$G$~contains no more than $T$ triangles. Then
$$
\alpha(G)
\ \ge \
\frac{n(G)}{10\Delta}
\left( \Strut{4ex}{0ex}
\log_2 \Delta \ - \ 
\sfrac{1}{2}\hspace{0.25ex}\log_2\!\left({\frac{T}{n(G)}}\right)\
\right)
$$
\end{lemma}

Observe that a graph has no triangles iff the neighborhood 
of every vertex is an independent set. If the neighborhood of
every vertex is sparse, then the graph will have few triangles.
This simple observation is made precise in the following theorem.
\vspace{-1.25ex}

\begin{theorem} 
\label{regular}
Let $G$ be a graph with maximum degree at most $\Delta$, 
and suppose that~for~all $v \in V(G)$,
the subgraph of $G$ induced by the neighborhood of\/ $v$ 
has at most\/ $t$ edges. Then
$$
\alpha(G)
\ \geq \
\frac{n(G)}{10\Delta}
\left( \Strut{3ex}{0ex}
\log_2 \Delta \ - \ 
\sfrac{1}{2}\hspace{0.25ex}\log_2\!\left({\frac{t}{3}}\right)\
\right)
$$
\end{theorem}

\Proof \looseness=-1
The number of triangles containing a given vertex $v \in V(G)$ is 
equal to the~number of edges in the subgraph of $G$\pagebreak[3.99]
induced by $N(v)$.
Thus, for every $v \in V(G)$, there are at most $t$ triangles containing $v$.
Summing the number of triangles containing $v$ over 
all $v \in V(G)$, we count each triangle in $G$ exactly three times. 
Hence, the total number of triangles in $G$ is at most $n(G)\,t/3$. 
The theorem now follows from \Lref{sparse}.
\qed\vspace{1ex}

Thus if we can show that $\GG$ 
is locally sparse
(that is, it satisfies~the~conditions of \Tref{regular} for
a relatively small value of~$t$), then we can improve upon the
\GV\ bound of \Pref{P-GV} by a factor of about 
$\log_2\!V(n,d{-}1)/10$.

\vspace{3.5ex}
\section{How Sparse Is the Sphere Graph?}
\label{sec4}

\noindent
In this section, 
we consider the Hamming sphere graph $\GS$, 
which is the subgraph of the Gilbert graph $\GG$ induced by the neighborhood
$N(\zero)$ of the vertex $\zero \in V(\GG)$. Clearly, the subgraph
induced in the Gilbert graph by the neighborhood $N(v)$ of any 
other vertex $v \in V(\GG)$ is isomorphic to $\GS$. Our goal
here is to determine how sparse $\GS$ is. Namely, we would like
to compute $e(\GS)$, the number of edges in $\GS$, and then 
determine~the~asymptotic relationship between $e(\GS)$ and 
the number of vertices in $\GS$. 
In view of \Lref{lemma-Gilbert} and \Tref{regular}, this would 
then provide a lower bound on $A_2(n,d) = \alpha(\GG)$.

For convenience, let us write $d' = d-1$. Recall that $\ceil{x}^+$ 
denotes the smallest {nonnegative} integer $m$ such that $m \ge x$,
for all $x \in \R$.
Consider the following simple lemma.\vspace{-1ex}
\begin{lemma}
\label{exact-degree}
Let $v \in V(\GS)$ be a vertex of weight $w$. Then the degree
of\/ $v$ in $\GS$ is given~by\vspace{-.75ex}
$$
\deg(v)
\ = \ 
\sum_{i=1}^{d'}
\sum_{j=\ceil{\frac{w+i-d'}{2}}^+}^{\min\{w,i\}}
{w \choose j} {n\,{-}\,w \choose i\,{-}\,j}
\ - \
1
\vspace{-.50ex}
$$
\end{lemma}

\Proof 
\looseness=-1
Let $u \in V(\GS)$ be a vertex of $\GS$ distinct from $v$, 
and suppose that $\wt(u) = i$
for some $i \in \{1,2,\ldots,d'\}$. Then 
$d(u,v) = \wt(u) + \wt(v) - 2|\chi(u) \cap \chi(v)|$,
where $\chi(\cdot)$ denotes the support of a vector in $\Fn$.
Write $j = |\chi(u) \cap \chi(v)|$. 
Then clearly $j \le \min\{w,i\}$.~~Fur\-thermore, 
$u$ and $v$ are adjacent in $\GS$ if and only
if $d(u,v) = w + i - 2j \le d'$. It follows that\vspace{-.75ex}
\be{exact-degree-aux}
\sum_{j=\ceil{\frac{w+i-d'}{2}}^+}^{\min\{w,i\}}
{w \choose j} {n\,{-}\,w \choose i\,{-}\,j}
\vspace{-.75ex}
\ee
is the number of vertices of weight $i$
that are adjacent to $v$, for all $i \ne w$.
For $i = w$, we need to subtract $1$ from the sum
in~\eq{exact-degree-aux}, because the sum counts
$v$ itself.
\qed\vspace{-1ex}

\begin{proposition}
\label{exact-e(GS)}\vspace{-1ex}
$$
e(\GS)
\ = \
\frac{1}{2}
\sum_{w=1}^{d'} \hspace{-.5ex} {n \choose w}\hspace{-.5ex}
\left(
\sum_{i=1}^{d'}\hspace{-.25ex}
\sum_{j=\ceil{\frac{w+i-d}{2}}^+}^{\min\{w,i\}}\hspace{-1ex}
{w \choose j}\! {n{-}w \choose i\,{-}\,j} - 1\hspace{-.5ex}
\right) 
$$
\end{proposition}

\Proof
Since $\GS$ has ${n \choose w}$ vertices of
weight $w$, this follows immediately from 
\Lref{exact-degree}.~\qed\vspace{1ex}

Comparing the foregoing expression for $e(\GS)$ 
with the expression for $e(n,d)$
in \Tref{exact}, we see that $e(n,d{-}1)$
is equal to $e(\GS)/3$. Thus \Pref{exact-e(GS)} in conjunction
with \Tref{regular} establish~\eq{th1}. 
This completes the proof of \Tref{exact}.\pagebreak[3.99]

\looseness=-1
Although \Pref{exact-e(GS)} gives an exact expression 
for $e(\GS)$, the asymptotic form of this expression is 
not immediately clear. Thus we now turn to
asymptotic bounds~on~$e(\GS)$. 
Observe that $|V(\GS)| = V(n,d') - 1$, so that
a complete graph on $V(\GS)$ has $\Omega\bigl(V(n,d')^2\bigr)$ edges. 
In contrast, 
we will show that under certain conditions,
there is an $\varepsilon > 0$ such that 
$e(\GS) = o\bigl(V(n,d')^{2-\varepsilon}\bigr)$.
To this end, the following lemma will be useful.
\vspace{-1ex}

\begin{lemma}
\label{increasing-degree}
Let $u$ and $v$ be vertices in $V(\GS)$
and suppose that $\wt(v) \leq \wt(u)$. 
Then $\deg(v) \geq \deg(u)$. 
\end{lemma} 
\Proof
\looseness=-1
It would suffice to prove that 
for all $w \in \{2,3,\ldots,d'\}$, we have $\deg(v) \geq \deg(u)$
if $\wt(u) = w$ and $\wt(v) = w-1$. 
Moreover, by \Lref{exact-degree} the degree of a vertex $u$ in $\GS$
depends on $u$ only through its Hamming weight $\wt(u)$. Thus we 
can assume without loss of generality that 
$$
\chi(u) = \{1,2,\ldots,w\}
\hspace{3ex}\text{and}\hspace{3ex}
\chi(v) = \{2,3,\ldots,w\}
$$
Now consider $N(u)$ and $N(v)$, the neighborhoods of $u$ and $v$
in~$\GS$. It is easy to see that\vspace{-.50ex}
\begin{eqnarray*}
N(u) - N(v) 
& \hspace{-1ex}=\hspace{-1ex} &
\{\,x \in V(\GS) ~:~ d(u,x) = d' \text{~and~} x_1 = 1 \,\}\\[.75ex]
N(v) - N(u) 
& \hspace{-1ex}=\hspace{-1ex} &
\{\,x \in V(\GS) ~:~ d(v,x) = d' \text{~and~} x_1 = 0 \,\}
\end{eqnarray*}
where $x_1$ denotes the first bit of the vector $x = (x_1,x_2,\ldots,x_n)$
in $V(\GS)$.
Let us denote the sets $N(u) - N(v)$ and $N(v) - N(u)$
by $\cA$ and $\cB$, respectively.
Let $\varphi: \Fn \to \Fn$ be the mapping 
$$
\varphi(x) \ = \ x + (100\cdots0)
$$
Note that $\varphi(u) = v$ and $\varphi(v) = u$. 
We claim that $\varphi(\cA) \subseteq \cB$.
Indeed, let us write $\varphi(x) = y = (y_1,y_2,\ldots,y_n)$. 
Evidently, if $d(u,x) = d'$ and $x_1 = 1$,
then $d(v,y) = d'$ and $y_1 = 0$. Moreover, for all $x \in \Fn$ with
$x_1 = 1$, the weight of $\varphi(x)$ is $\wt(x) - 1$. Thus if 
$x \in \cA$, then $\varphi(x) \in V(\GS)$ unless $x = (100\cdots0)$.
However $(100\cdots0) \,{\not\inn}\, \cA$, since the distance between 
$(100\cdots0)$ and $u$ is given by $w-1 \le d'-1 < d'$. This proves
that $\varphi(\cA) \subseteq \cB$. Since $\varphi$ is a bijection 
on $\Fn$, the fact that $\varphi(\cA) \subseteq \cB$ implies 
that $|\cA| \le |\cB|$. Hence $|N(v)| \ge |N(u)|$, and the lemma
follows. 
\qed\vspace{1.5ex}

The rest of our asymptotic analysis involves
the binary entropy function defined by
$$
H_2(x) \ \ \deff\ \ -x \log_2 \!x \,-\, (1-x) \log_2(1-x)
$$
for all $0 \le x \le 1$. In particular, we will make frequent use 
of the following lemma~\cite[pp.\,308-310]{MWS}, 
which is a well-known estimate for a sum of binomial
coefficients.\vspace{-1ex}
\begin{lemma}
\label{binom-H2}
Let $\mu \in \R$, 
and suppose that $\mu n$ is an integer
in the range $1 \le \mu n \le 0.5n$. Then
\be{eq:binom-H2}
\frac{2^{nH_2(\mu)}}{\sqrt{8n\mu(1{-}\mu)}}
\ \le \
\sum_{k = 0}^{\mu n} {n \choose k} 
\ \le \
2^{nH_2(\mu)}
\ee
\end{lemma}

Now, let $\lambda$ be a real number in the 
range $\sfrac{2}{3}\le \lambda < 1$.
To~simplify notation,
we~henceforth assume that $d' < 0.5n$ and that~$\lambda d'$ is an integer
(this obviates the need for numerous 
$\ceil{\cdot}$ and $\trunc{\cdot}$
functions in what follows).
We will derive a bound on $e(\GS)$ by considering 
separately vertices of weight $< \lambda d'$ and 
vertices of weight $\ge \lambda d'$ in $\GS$. Thus we write
\be{split}
e(\GS) \ + \ \frac{1}{2} \sum_{i=1}^{d'} {n \choose i}
\ = \ 
\frac{e_1(\lambda,n,d) \ + \ e_2(\lambda,n,d)}{2}
\vspace{-.5ex}
\ee
with
\begin{eqnarray}
\hspace*{-1.50ex}e_1(\lambda,n,d)
& \hspace{-1.0ex}{\deff}\hspace{-1ex} &
\hspace{-1.5ex} \sum_{v \inn V(\GS) \atop \wt(v) < \lambda d'}
\hspace{-1ex} (\deg(v) + 1)
\nonumber\\
\label{e1}
& \hspace{-2ex}{=}\hspace{-1ex} &
\hspace{-1.0ex}
\sum_{w=1}^{\lambda d' - 1}
\hspace{-.50ex}
{n \choose w}
\sum_{i=1}^{d'}
\sum_{j=\ceil{\frac{w+i-d'}{2}}^+}^{\min\{w,i\}}
\hspace{-1ex}{w \choose j} {n{-}w \choose i{-}j}
\\[.5ex]
\hspace*{-1.50ex}e_2(\lambda,n,d)
& \hspace{-1.5ex}{\deff}\hspace{-1ex} &
\hspace{-1.5ex} \sum_{v \inn V(\GS) \atop \wt(v) \ge \lambda d'}
\hspace{-1ex} (\deg(v) + 1)
\nonumber\\
\label{e2}
& \hspace{-2ex}{=}\hspace{-1ex} &
\hspace{-1.50ex}
\sum_{w=\lambda d'}^{d'}
\hspace{-.50ex}
{n \choose w}
\sum_{i=1}^{d'}
\sum_{j=\ceil{\frac{w+i-d'}{2}}^+}^{\min\{w,i\}}
\hspace{-1ex}{w \choose j} {n{-}w \choose i{-}j}\hspace*{3ex}
\end{eqnarray}
where the explicit expressions for $e_1(\lambda,n,d)$ and
$e_2(\lambda,n,d)$ 
follow from \Lref{exact-degree} and \Pref{exact-e(GS)}.
Let $v$ be a vertex in $V(\GS)$ with $\wt(v) = 1$. 
Note that 
$\ceil{(1+i-d')/2}^+ \! = 0$~for~all $i \le d'\,{-}\,1$.
Therefore, by \Lref{exact-degree}, we have
\begin{eqnarray*}
\deg(v) + 1
& \hspace{-1ex}{=}\hspace{-1ex} &
\sum_{i=1}^{d'-1} 
\left(\Strut{4.5ex}{0.10cm} {n{-}1 \choose i} + {n{-}1 \choose i{-}1} \right)
\ + \ {n{-}1 \choose d'{-}1}
\\
& \hspace{-2ex}{\le}\hspace{-1ex} &
\sum_{i=1}^{d'} {n \choose i} 
\ \le \
2^{n H_2(\delta)}
\end{eqnarray*}
where $\delta = d'/n$ and the last inequality
follows from \Lref{binom-H2}. Combining the 
definition of $e_1(\lambda,n,d)$ in~\eq{e1} with \Lref{increasing-degree}
thus produces the following bound
\be{e1-bound}
e_1(\lambda,n,d)
\ \le \
2^{n H_2(\delta)} \sum_{w=1}^{\lambda d' - 1} {n \choose w}
\ \le \
2^{n \bigl(H_2(\delta) + H_2(\lambda\delta)\bigl)}
\ee
Now, let $v$ be a vertex in $V(\GS)$ with $\wt(v) = \lambda d'$.
Then, again by \Lref{exact-degree}, the degree of $v$ in $\GS$
is given by 
$$
\deg(v) + 1
\ = \
h_1(\lambda,n,d) \,+\, h_2(\lambda,n,d) 
\vspace{-.5ex}
$$
with\vspace{-2ex}
\begin{eqnarray}
\label{f1}
h_1(\lambda,n,d)
& \hspace{-1.0ex}{\deff}\hspace{-1ex} &
\sum_{i=1}^{\mu d'} \sum_{j=0}^i
{w \choose j} \hspace{-0.50ex} {n{-}w \choose i - j}
\\[1ex]
\label{f2}
h_2(\lambda,n,d)
& \hspace{-1.0ex}{\deff}\hspace{-1ex} &
\hspace{-1ex}\sum_{i=\mu d'+1}^{w - 1}
\hspace{-.25ex}\sum_{j=\ceil{\frac{i - \mu d'}{2}}}^i
\hspace{-1ex}{w \choose j} \hspace{-0.50ex} {n{-}w \choose i - j}
\nonumber\\
& & \ + \
\sum_{i=w}^{d'}
\hspace{-.25ex}\sum_{j=\ceil{\frac{i - \mu d'}{2}}}^{w}
\hspace{-1ex}{w \choose j} \hspace{-0.50ex} {n{-}w \choose i - j}
\end{eqnarray}
where we have introduced the notation
\mbox{$w = \lambda d'$ and $\mu = 1 - \lambda$}.
To upper-bound $h_1(\lambda,n,d)$, 
observe that for all $i$ and $j$\pagebreak[3.99]
in the double-sum of~\eq{f1}, we have $j \le i \le \mu d' \le 0.5 w$
and therefore
${w \choose j} \le {w \choose i}$. Thus\vspace{-1.50ex}
\begin{eqnarray}
h_1(\lambda,n,d)
& \hspace{-1.0ex}{\le}\hspace{-1ex} &
\sum_{i=1}^{\mu d'} {w \choose i}
\sum_{j=0}^i {n{-}w \choose i - j}
\nonumber\\[.05ex]
& \hspace{-1.0ex}{\le}\hspace{-1ex} &
\sum_{i=0}^{\mu d'} {w \choose i}
\sum_{j=0}^{\mu d'} {n{-}w \choose j}
\nonumber\\[.75ex]
\label{f1-bound}
& \hspace{-1.0ex}{\le}\hspace{-1ex} &
2^{n \lambda \delta H_2\left({\mu\over\lambda}\right) \ + \
n(1-\lambda\delta) H_2\left({\mu\delta \over 1 - \lambda\delta}\right)}
\end{eqnarray}
To upper-bound $h_2(\lambda,n,d)$ in~\eq{f2}, we will use the 
trivial estimate 
${w \choose j} \le 2^{w} = 2^{n \lambda\delta}$ for all $j$
(in the case of~\eq{f2}, this estimate is actually not too far off).
Thus
\begin{eqnarray}
h_2(\lambda,n,d)
& \hspace{-1.0ex}{\le}\hspace{-1ex} &
2^{n \lambda\delta} \hspace{-1ex}
\sum_{i=\mu d'+1}^{w - 1}
\hspace{-.25ex}\sum_{j=\ceil{\frac{i - \mu d'}{2}}}^i
\hspace{-1ex} {n{-}w \choose i - j}
\nonumber\\[.25ex]
\label{f2-aux}
& & 
 + \ \
2^{n \lambda\delta} 
\sum_{i=w}^{d'}
\hspace{-.25ex}\sum_{j=\ceil{\frac{i - \mu d'}{2}}}^{w}
\hspace{-1ex} {n{-}w \choose i - j}
\end{eqnarray}
Since the summation on $j$ in the second double-sum of~\eq{f2-aux}
is up to $w \le i$, we can proceed with the upper bound by uniting
the two double-sums as follows
\begin{eqnarray}
h_2(\lambda,n,d)
& \hspace{-1.0ex}{\le}\hspace{-1ex} &
2^{n \lambda\delta} \hspace{-1.5ex}
\sum_{i=\mu d'+1}^{d'}
\hspace{-.25ex}\sum_{j=\ceil{\frac{i - \mu d'}{2}}}^i
\hspace{-1ex} {n{-}w \choose i - j}
\nonumber\\[.25ex]
\label{f2-aux3}
& \hspace{-1.0ex}{=}\hspace{-1ex} &
2^{n \lambda\delta} \hspace{-1.5ex}
\sum_{i=\mu d'+1}^{d'}
\hspace{-.75ex}\sum_{j=0}^{\trunc{\frac{i + \mu d'}{2}}}
\hspace{-1ex} {n{-}w \choose j}
\end{eqnarray}
where the equality in~\eq{f2-aux3} follows by a straightforward
change of variables.
Finally, observing that  
$(i + \mu d')/2 \le d' - \frac{\lambda}{2}d'$
for all $i \le d'$, we get\vspace{-1ex}
\begin{eqnarray}
h_2(\lambda,n,d)
& \hspace{-1.0ex}{\le}\hspace{-1ex} &
2^{n \lambda\delta} \hspace{-1.5ex}
\sum_{i=\mu d'+1}^{d'}
\hspace{-1.00ex}\sum_{j=0}^{\trunc{d' - \frac{\lambda}{2}d'}}
\hspace{-1ex} {n {-} w \choose j}
\nonumber\\[.75ex]
\label{f2-bound}
& \hspace{-1.0ex}{\le}\hspace{-1ex} &
{n \lambda\delta}\
2^{n \lambda\delta \ + \ n(1-\lambda\delta) 
H_2\left({\delta - \frac{\lambda}{2}\delta \over 1-\lambda\delta} \right)}
\end{eqnarray}
Combining~\eq{f1-bound} and~\eq{f2-bound} with the definition 
of $e_2(\lambda,n,d)$ in~\eq{e2} and once again invoking
\Lref{increasing-degree}, we obtain the following bound\vspace{-1.25ex}
\begin{eqnarray}
\hspace*{-4ex}
e_2(\lambda,n,d)\hspace*{-.5ex}
& \hspace{-1.0ex}{\le}\hspace{-1.25ex} &
\hspace{-.25ex}\bigl(\deg(v) + 1\bigr) \sum_{w=\lambda d'}^{d'} {n \choose w}
\\[.25ex]
& \hspace{-1.0ex}{\le}\hspace{-1.25ex} &
\hspace{-.5ex}\Bigl(h_1(\lambda,n,d) + h_2(\lambda,n,d)\Bigr) 
\sum_{w=0}^{d'} {n \choose w}
\\[.75ex]
& \hspace{-1.0ex}{\le}\hspace{-1.25ex} &
\label{e2-aux}
2^{n 
\biggl(
H_2(\delta) + 
\lambda \delta H_2\left({\mu\over\lambda}\right) + 
(1-\lambda\delta) H_2\left({\mu\delta \over 1 - \lambda\delta}\right)
\biggr)}
\nonumber\\ 
& & + \:\
{n \lambda\delta}\,
2^{n 
\biggl(
H_2(\delta) + 
\lambda\delta + 
(1-\lambda\delta) 
H_2\left({\delta - \lambda\delta/2 \over 1-\lambda\delta} \right)
\biggr)}
\\[.5ex]
\label{e2-bound}
& \hspace{-1.0ex}{\le}\hspace{-1.25ex} &
\hspace*{-0.50ex}(n \lambda\delta \,{+}\, 1)\,
2^{\!n 
\biggl(\hspace{-0.5ex}
H_2(\delta) + 
\lambda\delta + 
(1-\lambda\delta) 
H_2\left({\delta - \lambda\delta/2 \over 1-\lambda\delta} \right)
\hspace{-0.5ex}\biggr)}
\end{eqnarray} 
where \eq{e2-bound} follows from the fact that for 
\mbox{$\,\sfrac{2}{3} \,{<}\,\lambda\,{<}\, 1$}~and~\mbox{$\,\delta \le 0.5$},
the first exponent in~\eq{e2-aux} is
strictly less than the second exponent. We are now
ready to prove the following proposition.
\vspace{-1ex}

\begin{proposition}
\label{P-asymptotic}
Let $\varepsilon$ and $\lambda$ be positive real numbers
strictly less than $1$, with $\lambda \ge \sfrac{2}{3}$. 
Then $e(\GS) = o\bigl(V(n,d')^{2-\varepsilon}\bigr)$,
providing $\delta = d'/n$ satisfies the following two conditions:
\begin{eqnarray}
\label{condition1}
(1 {-} \varepsilon)\, H_2(\delta) 
& \hspace{-1.0ex}{>}\hspace{-1.00ex} &
H_2(\lambda \delta) 
\\[.25ex]
\label{condition2}
(1 {-} \varepsilon)\, H_2(\delta) 
& \hspace{-1.0ex}{>}\hspace{-1.00ex} &
\lambda \delta \, + \, 
(1{-}\lambda\delta)\,
H_2\!\left({\delta - \lambda\delta/2 \over 1-\lambda\delta} \right)
\end{eqnarray} 
\end{proposition}
\Proof
We estimate $e(\GS)$ by combining~\eq{split} with the upper
bounds in \eq{e1-bound}~and~\eq{e2-bound}~on $e_1(\lambda,n,d)$
and $e_2(\lambda,n,d)$. It follows that the ratio 
$e(\GS)/V(n,d')^{2-\varepsilon}$
is upper-bounded by
\begin{eqnarray*}
\lefteqn{\frac{e(\GS)}{V(n,d')^{2-\varepsilon}}
\,\ \le \,\
\frac{2^{n \bigl(H_2(\lambda\delta) \ - \ (1-\varepsilon)H_2(\delta)\bigl)}}
{\left(\Strut{2ex}{0.1ex}8n\delta(1{-}\delta)\right)^{{\varepsilon\over2}-1}}}
\\[1.00ex]
& &  + \:\
\frac{(n \lambda\delta + 1)\,
2^{n 
\left(
\lambda\delta \ + \ 
(1-\lambda\delta) 
H_2\left({\delta - \lambda\delta/2 \over 1-\lambda\delta}\right)
\ - \ 
(1-\varepsilon)H_2(\delta) 
\Strut{2ex}{0.1ex}\right)}}
{\left(\Strut{2ex}{0.1ex}8n\delta(1{-}\delta)\right)^{{\varepsilon\over2}-1}}
\end{eqnarray*}
where we used \Lref{binom-H2} to bound $V(n,d')$.
It is clear that if
$\delta$ satisfies \eq{condition1}\,--\,\eq{condition2},
then the right-hand side of the above expression tends to zero
(exponentially fast) as $n \to \infty$.~\qed\vspace{1ex}

Motivated by \Pref{P-asymptotic}, we now introduce the functions
$f_{\varepsilon,\lambda}(\delta)$ and $g_{\varepsilon,\lambda}(\delta)$
with domain $\delta \in [0,0.5]$,
parametrized by $\varepsilon$ and $\lambda$ and defined as follows
\begin{eqnarray*}
f_{\varepsilon,\lambda}(\delta)
& \hspace{-1.0ex}{\deff}\hspace{-1.0ex} &
(1 {-} \varepsilon)\, H_2(\delta) 
\ - \
H_2(\lambda \delta) 
\\
g_{\varepsilon,\lambda}(\delta)
& \hspace{-1.0ex}{\deff}\hspace{-1.0ex} &
(1 {-} \varepsilon)\, H_2(\delta) 
\ - \
\lambda \delta \ - \ 
(1{-}\lambda\delta)\,
H_2\!\left({\delta - \lambda\delta/2 \over 1-\lambda\delta} \right)
\end{eqnarray*}
The two functions $f_{\varepsilon,\lambda}(\delta)$ and 
$g_{\varepsilon,\lambda}(\delta)$ are plotted in Figure\,1
and Figure\,2, respectively,~for $\varepsilon = 0.000001$
and $\lambda = 0.999$. Figure\,3 shows a close-up view of
these two functions (for the same $\varepsilon$ and $\lambda$)
in the range $\delta \in [0.499,0.5]$. It can be seen from 
Figures 1, 2, and 3 that conditions \eq{condition1} and
\eq{condition2} of \Pref{P-asymptotic} are satisfied whenever
$d'/n \le 0.4994$.

\looseness=-1
We are now ready to complete the proof of \Tref{asymptotic}.
By \Lref{lemma-Gilbert}, $A_2(n,d) = \alpha(\GG)$,
where $\GG$ is the Gilbert graph defined in Section\,\ref{sec3}.
The Gilbert graph is a $\Delta$-regular graph on $|V(\GG)| = 2^n$ 
vertices with constant degree $\Delta = V(n,d') - 1$.
The subgraph of~$\GG$ induced by the neighborhood of any vertex in $V(\GG)$
is isomorphic to the sphere graph~$\GS$ and has exactly $e(\GS)$
edges. Therefore, by \Tref{regular}, for all $\varepsilon > 0$, 
we have\vspace{.75ex}
\begin{eqnarray*}
A_2(n,d)
& \hspace{-1.0ex}{\ge}\hspace{-1.0ex} &
\frac{2^n}{V(n,d')} 
\cdot
\frac{\log_2\! V(n,d') - \sfrac{1}{2}\log_2 e(\GS)}{10}
\\[2ex]
& \hspace{-1.0ex}{=}\hspace{-1.0ex} &
\frac{2^n}{V(n,d')} 
\left(
\frac{\varepsilon\log_2\! V(n,d')}{20}
\right.
\nonumber\\
\label{th2-proof}
& & \hspace*{12ex} + \
\left.
\frac{\log_2\! V(n,d')^{2{-}\varepsilon} - \log_2\! e(\GS)}{20}
\right)
\nonumber
\end{eqnarray*}
By \Pref{P-asymptotic}, the ratio $e(\GS)/V(n,d')^{2-\varepsilon}$
tends to zero for $\varepsilon = 0.000001$, whenever $d'/n < d/n \le 0.4994$
(cf.~Figure\,3).\pagebreak[3.99]
Therefore, the second fraction in parentheses 
becomes positive for all sufficiently large $n$,
and \Tref{asymptotic} follows.\pagebreak[4.00]

\centerline{\psfig{figure=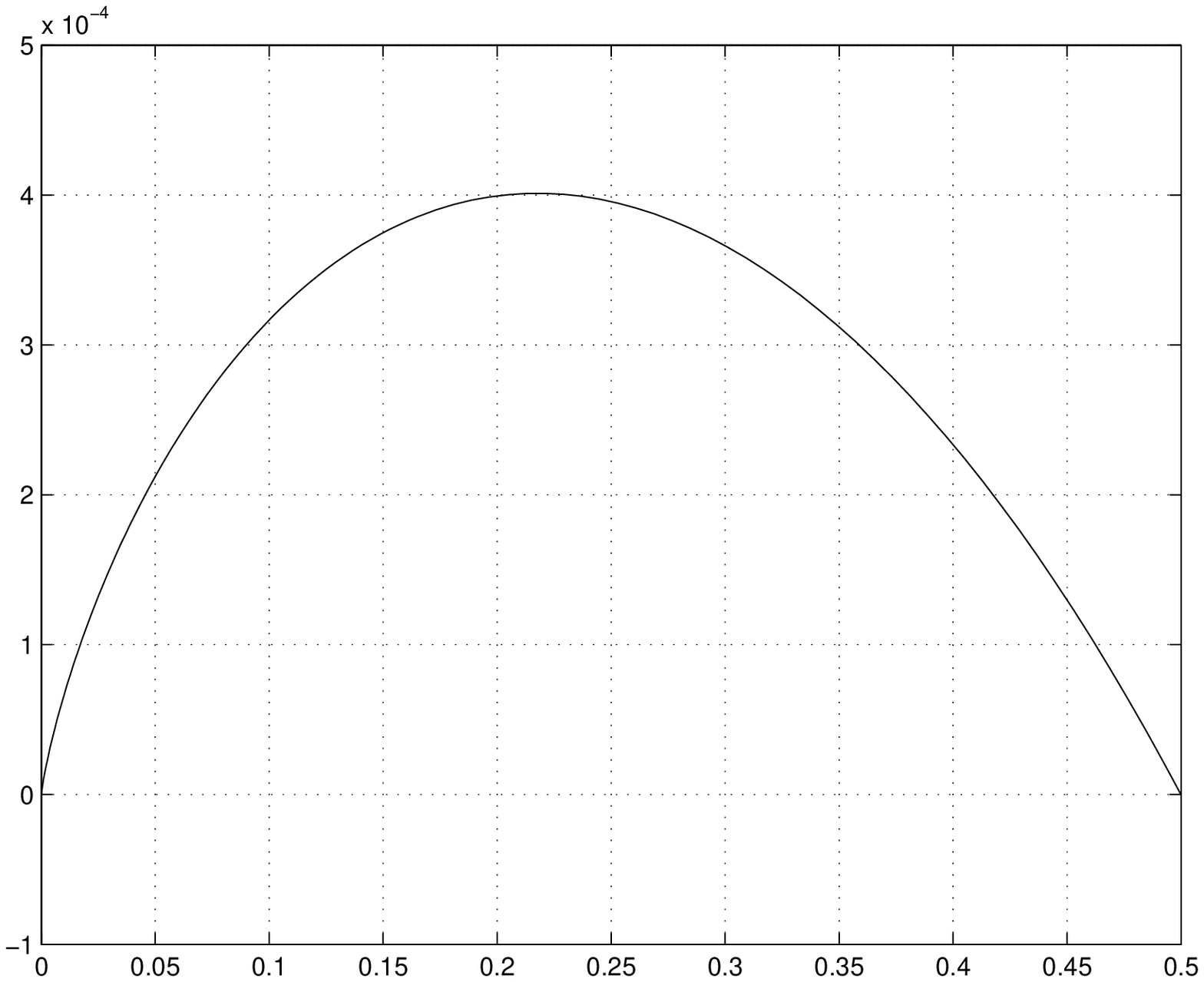,width=2.50in,silent=}}

\begin{center}
{\small\bf Figure\,1.}  
{\small\sl Plot of the function $f_{\varepsilon,\lambda}(\delta)$
for $\varepsilon = 0.000001$ and $\lambda = 0.999$}
\vspace{3.0ex}
\end{center}

\centerline{\psfig{figure=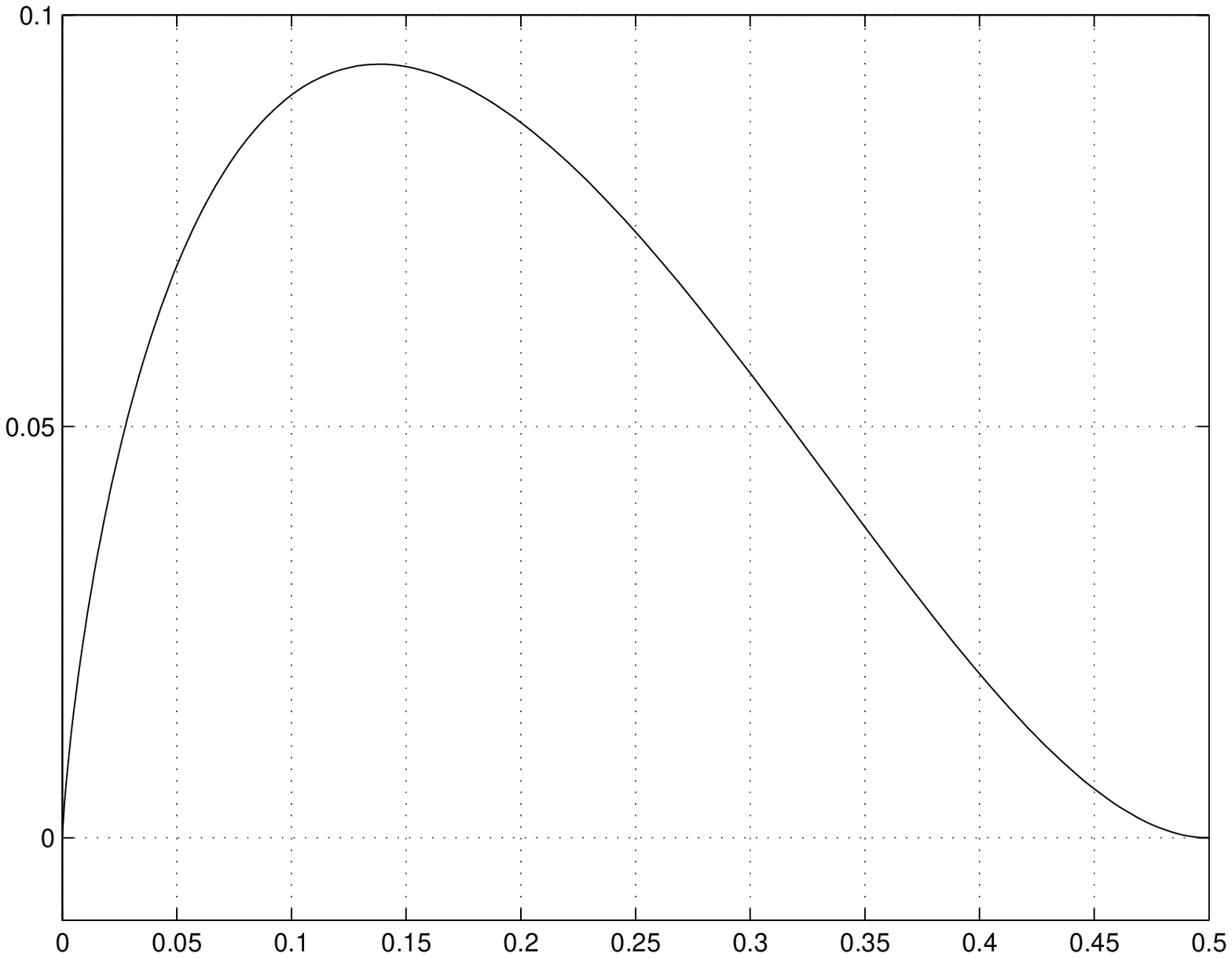,width=2.50in,silent=}}

\begin{center}
{\small\bf Figure\,2.}  
{\small\sl Plot of the function $g_{\varepsilon,\lambda}(\delta)$
for $\varepsilon = 0.000001$ and $\lambda = 0.999$}
\vspace{3.0ex}
\end{center}

\centerline{\psfig{figure=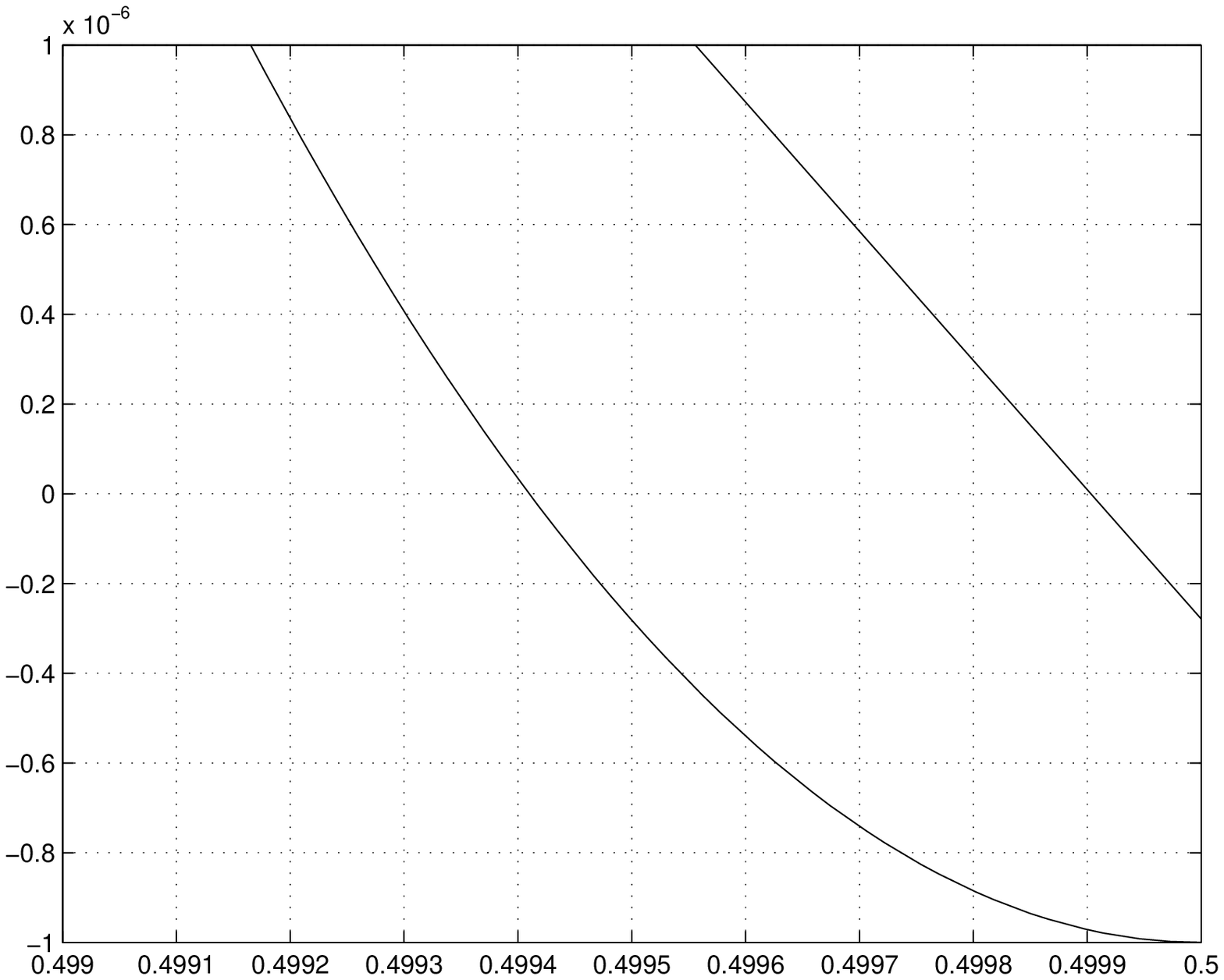,width=2.50in,silent=}}

\begin{center}
{\small\bf Figure\,3.} \vspace{-.5ex} 
{\small\sl Close-up view of the functions $f_{\varepsilon,\lambda}(\delta)$
and $g_{\varepsilon,\lambda}(\delta)$ in the neighborhood of $\delta = 0.5$
for $\varepsilon = 0.000001$ and $\lambda = 0.999$}
\vspace{3.0ex}
\end{center}

\noindent\looseness=-1
{\bf Remark.}
We note that the degree of a vertex $v$ in $\GS$ is related
to the so-called \emph{intersection numbers} $p^w_{i,k}$ of the Hamming
association scheme $\cH(n,2)$ --- 
see~\cite{Delsarte} and \cite[Chapter\,21]{MWS}
for a~detailed description of $\cH(n,2)$. 
Specifically, given any two vectors $u,v \in \Fn$ with $d(u,v) = w$,
the intersection number $p^w_{i,k}$ is defined as the number 
of vectors $x \in \Fn$ such that $d(x,u) = i$ and $d(x,v) = k$.
Thus the sum \eq{exact-degree-aux} can be written as
$$
p^w_{i,0} + p^w_{i,1} + \cdots + p^w_{i,d'}
$$ 
An explicit expression for $p^w_{i,k}$
is given in~\cite[p.\,656]{MWS}.
However, the proof of \Lref{exact-degree} above, 
which does not use the intersection numbers, appears to
be simpler and shorter.\vspace{1ex}


\noindent\looseness=-1
{\bf Remark}. 
To get the best threshold on $d/n$ such that \eq{th2} holds,
one should optimize the value of $\lambda$ for a given $\varepsilon$
in \Pref{P-asymptotic}
(alternatively, one could try to directly find the maximum 
term in the triple-sum of \Pref{exact-e(GS)}). 
We have made no special effort to optimize
this threshold beyond $0.499$. However,~we believe that with an
appropriate choice of $\varepsilon$, $\lambda$
in \Pref{P-asymptotic} (or with other methods),
one can get as close as desired to the ultimate threshold $d/n \le 0.5$.
It is, therefore, 
surprising that for $\delta = 0.5$, the
number of edges in $\GS$ is very close to $V(n,d')^2$.\vspace{-.5ex}
\begin{proposition}
\label{0.5-theta}
If $d'/n = 0.5$, then $e(\GS) \ge 0.25 V(n,d')^2$.
\vspace{-.5ex}
\end{proposition}

\Proof
Let $v \in V(\GS)$ be a vertex of weight $d' = n/2$,
and let $\one$ denote the all-one vector $(11\cdots1)$ in $\Fn$.
Then $\one + v$ is another vertex of weight $d'$ in $V(\GS)$.
Given any other vertex $u \in V(\GS)$, we have 
$d(u,v) + d(u,\one{+}v) = n = 2d'$, so that $u$ is 
adjacent to at least one of $v$ or $\one + v$.
Thus every vertex in $\GS$ is adjacent to at least half 
of the vertices of weight $d'$ (excluding, possibly, itself). This
implies that\vspace{-.75ex}
\be{remark-aux}
\sum_{v \inn V(\GS) \atop \wt(v) = d'} \!\!
(\deg(v) + 1)
\ \ge \
\frac{1}{2} 
\sum_{w=1}^{d'} {n \choose w} {n \choose d'}
\vspace{-.75ex}
\ee
By \Lref{exact-degree}, all the vertices of weight $d'$ have
the same degree in $\GS$. Thherefore, it follows from~\eq{remark-aux}
that for {\sl every\/} $v \in V(\GS)$
with $\wt(v) = d'$, we have\vspace{-1.50ex}
\be{deg-v-heavy}
\deg(v) + 1 
\ \ge \
\frac{1}{2} \sum_{w=1}^{d'} {n \choose w}
\vspace{-.50ex}
\ee
Now, by \Lref{increasing-degree}, the degree of all other
vertices in $\GS$ is greater or equal to the degree of a vertex
of weight $d'$. This essentially establishes the proposition.
It remains to worry about the fact that the sum on the right-hand
side of~\eq{deg-v-heavy} does not include the term ${n \choose 0}$
and about the extra $1$ on the left-hand side of~\eq{deg-v-heavy}.
We omit these tedious details.
\qed\vspace{1ex}

\noindent\looseness=-1
Thus it appears that the sphere graph $\GS$ transitions abruptly
from being sparse to being nearly complete at $d'/n = 0.5$.
We do not have an intuitive ``explanation'' for this phenomenon,
but note that it is reminiscent of threshold phenomena for 
codes and graphs observed in~\cite{ZC} and in~\cite{Margulis},
respectively.\vspace{1.00ex}

We also note 
that for $d/n \ge 0.5$, 
the problem of determining $\And$ is essentially settled.
Provided enough Hadamard matrices exist, 
$A_2(2d,d) = 4d$ and 
$A_2(n,d) = 2\trunc{d/(2d-n)}$
for all even $d$ with $2d \,{>}\, n$.
This is the well-known result of
Levenshtein~\cite{Levenshtein61},
who constructed codes achieving the Plotkin 
bound~\cite[pp.\,41--43]{MWS} from Hadamard 
matrices.

\vspace{3.00ex}
\section{Generalizations and Open Problems}
\label{sec5}

\noindent\looseness=-1
The well-known proofs by Gilbert~\cite{Gilbert} and
Varshamov~\cite{Varshamov} of the bounds 
in~\eq{2-GV} and~\eq{Varshamov-f}, respectively, are 
``constructive'' in that they provide simple (but exponential-time) 
algorithms to construct codes whose parameters meet or 
exceed the corresponding bounds.
Moreover, Gilbert's ``constructive'' argument~\cite{Gilbert} has 
been extended to quite general contexts~\cite{GF,Tolhuizen,Ytrehus}
using the so-called altruistic algorithm
(which is also exponential-time). 

\looseness=-1
We would like to point out that the bound 
of \Tref{asymptotic} is ``constructive'' in the
same sense as~\cite{GF,Tolhuizen,Varshamov}, 
and~\cite{Ytrehus}.\pagebreak[3.99]
Hofmeister and Lefmann~\cite{HL} provide
an algorithm which, given any $\Delta$-regular graph $G$ 
with at most $n(G)\Delta^{2-\varepsilon}$
triangles, finds an independent set of size 
at least 
$\Omega(n(G) \log_2(\Delta)/\Delta)$ in~$G$.
By the results of Section\,IV,  
the Gilbert graph $\GG$ contains at most $O(n(\GG) \Delta^{2-\varepsilon})$
triangles whenever $d/n \le 0.499$.
Thus, when applied to the Gilbert graph $\GG$, 
the Hofmeister-Lefmann algorithm~\cite{HL} will
produce codes satisfying~\eq{th2}. 
The Hofmeister-Lefmann algorithm runs in time that is polynomial 
in the size of $\GG$ but, of course, exponential in 
the code length $n$.\vspace{1.5ex}

Up to now, for the sake of brevity, we have focused exclusively
on binary codes. Nevertheless, it should be clear that Theorems
\ref{exact} and \ref{asymptotic} can be generalized to arbitrary
alphabets of size $q$, where $q$ need not even be a prime power.
Here, we give a generalization to $q$-ary alphabets 
of \Tref{exact}.\vspace{-1ex}
\begin{theorem}
\label{q-ary}
Let $q$, $n$, and $d$ be positive integers 
with $d \le n$ and $q \ge 2$.
Define the volume~of a $q$-ary Hamming sphere of radius~$d$ as
$V_q(n,d) = \sum_{i=0}^d {n \choose i} (q{-}1)^i$, and let\hspace*{1.20in}
\begin{eqnarray*}
e_q(n,d) 
& \hspace{-2.00ex}\deff\hspace{-2.00ex} &
\frac{1}{6}\hspace{-.28ex}
\sum_{w=1}^{d}\hspace{-.28ex}
\sum_{i=1}^{d}\hspace{-.28ex}
\sum_{j=1}^{a}\hspace{-.28ex}
\sum_{k=b}^{a-j}\hspace{-.25ex}
\frac{n!\,(q{-}2)^k(q{-}1)^{w+i-c}}%
{j!\hspace{.3ex}
k!\hspace{.3ex}
(w{-}c)!\hspace{.3ex}
(i{-}c)!\hspace{.3ex}
(n{+}c{-}w{-}i)!}
\\
& & \hspace*{5ex}- \
\frac{V_q(n,d)}{6} 
\end{eqnarray*}
where $a \ \deff\, \min\{w,i\}$, $c \ \deff\, j+k$, and $b$ is 
the smallest {nonnegative} integer that is greater or equal to 
$(w{+}i) - j - \min\{d{+}j,n\}$.
Then\vspace{-1ex} 
\be{th14}
A_q(n,d)\hspace{-.1ex} 
\ \ge \
\hspace{-.1ex}\frac{q^n}{V_q(n,d{-}1)} 
\cdot
\frac{\log_2\! V_q(n,d{-}1) \hspace{.2ex}-\hspace{.2ex} 
\log_2\! \sqrt{e_q(n,d{-}1)}}{10}\vspace{-.25ex}
\ee
\end{theorem}

\Proof
Let $\cA$ be an alphabet with $q$ letters. We define 
the $q$-ary Gilbert graph $\GGq$ as before, namely
$V(\GGq) = \cA^n$ 
and 
$\{u,v\} \in E(\GGq)$ if and only if\, $1 \le d(u,v) \le d'$.
Then $\GGq$ is $\Delta$-regular with $\Delta = V_q(n,d') - 1$,
and \Tref{regular} applies. It remains to count the number 
of edges in the graphs induced in $\GGq$ by the neighborhoods
of its vertices. 
Without loss of generality,
we can call any one of the letters of $\cA$ ``zero,'' 
and consider the graph $\GSq$ which is induced in the  
$q$-ary Gilbert graph by the neighborhood $N(\zero)$ 
of the vertex $\zero \in \cA^n$.

\looseness=-1
Let $u = (u_1,u_2,\ldots,u_n)$ and 
$v = (v_1,v_2,\ldots,v_n)$~be~two~ver\-tices 
of $\GSq$ with $\wt(v) = w$ and $\wt(u) = i$
(observe that Ham\-ming weight is well-defined, once we have
identified a~$0 \in \cA$). Let\vspace{-1.5ex}
\begin{eqnarray*}
j 
& \hspace{-.75ex}\deff\hspace{-.75ex} &
\bigl|\{l : u_l = v_l, v_l \ne 0, u_l \ne 0\}\bigr|\\[0.50ex]
k 
& \hspace{-.75ex}\deff\hspace{-.75ex} &
\bigl|\{l : u_l \ne v_l, v_l \ne 0, u_l \ne 0\}\bigr|
\end{eqnarray*}
It is obvious that $j \le \min\{w,i\}$.
Further, if $j$ is already fixed, then clearly $k \le \min\{w,i\} - j$.
Moreover, $w - k - j$ is the number of positions $l$ such that
$v_l \ne 0$ and $u_l = 0$. This number cannot be greater than
$n - \wt(u) = n - i$, which implies that 
$
k \ge (w{+}i) - j - n
$.
Finally, it is easy to see that 
$$
d(u,v) \ = \ w+i-2j-k
$$
so that the vertices 
$u$ and $v$ are adjacent in $\GSq$ if and only if
$
k \ge (w{+}i) - j - (d'{+}j)
$.
Putting all this together,~we~can~enume\-rate
the total number of vertices of weight 
\mbox{$i \ne w$} that~are adja- cent in $\GSq$ 
to a fixed vertex $v \in V(\GSq)$
of weight $\wt(v) = w$ as follows
\be{degv-i-qary}
\sum_{j=0}^a
\sum_{k=b}^{a-j}
{w \choose j} {w{-}j \choose k} (q{-}2)^k 
{n{-}w \choose i-c} (q{-}1)^{i-c}
\ee
where $a$, $b$, and $c$ are as defined in the theorem.
For $i=w$, we again need to subtract $1$ from the sum
in~\eq{degv-i-qary} since the sum then counts $v$ itself.
Enumerating over all possible values of $i$, we find 
that the degree of $v$ in $\GSq$ is given  by
\be{degv-qary}
\sum_{i=1}^{d'}
\sum_{j=0}^a
\sum_{k=b}^{a-j}
{w \choose j} {w{-}j \choose k} {n{-}w \choose i-c} 
(q{-}2)^k (q{-}1)^{i-c}
\ - \ 
1
\ee
The total number of vertices of weight $w$ in $\GSq$ 
is ${n \choose w} (q{-}1)^w$. Combining this with~\eq{degv-qary}
produces an expression for $e(\GSq)$, and it is easy
to see that $e_q(n,d{-}1) = e(\GSq)/3$.
\qed\vspace{1.5ex}

\noindent
{\bf Remark.}
We could have used the {intersection numbers~$p^w_{i,j}$}
of the Hamming association scheme $\cH(n,q)$ in the proof of
\Tref{q-ary}. Specifically, the sum in \eq{degv-i-qary} can 
again be written as
$
p^w_{i,0} + p^w_{i,1} + \cdots + p^w_{i,d'}
$. 
Therefore 
$$
e(\GSq)
\ = \
\frac{1}{2}
\sum_{w=1}^{d'}
\sum_{i=1}^{d'}
\sum_{j=1}^{d'} 
{n \choose w}\, p^w_{i,j}\, (q{-}1)^w
$$
with the convention that $p^w_{i,j} = 0$ when $w > i+j$.
A formula for the intersection numbers of the $q$-ary Hamming scheme
$\cH(n,q)$ may be found 
in~\cite[eq.\,(2)]{BGS}. While the resulting expression for
$e_q(n,d)$ is shorter than its counterpart in \Tref{q-ary},
we prefer the latter since it is more explicit.
\vspace{1.00ex}

In the original version of this paper, we have left
the asymptotic investigation of the bound in \Tref{q-ary} 
as an open problem, and conjectured that it should
lead to 
\be{q-conjecture}
A_q(n,d) 
\ \ge\
c\, \frac{q^n}{V_q(n,d{-}1)}\, \log_2\! V_q(n,d{-}1)
\ee
for some positive constant $c$. 
This conjecture has been~proved in the recent
work of Vu and Wu~\cite{VuWu}. Specifically,
Vu and Wu~\cite{VuWu} show that if
$d/n < (q{-}1)/q$, then~\eq{q-conjecture} 
holds for a constant $c$ that depends on the ratio $d/n$.
They also give an explicit, though rather elaborate,
expression for $c$ in terms of $d/n$.
\vspace{1ex}

Our general approach can be extended to 
many other situations where generalizations of
the Gilbert-Varshamov bound are currently used. We give
just one concrete example.

Let $\Andw$ denote the the maximum number of codewords in
a binary code of length~$n$, constant weight $w$, and
minimum Hamming distance $2d$. Levenshtein~\cite{Levenshtein71} has 
generalized the Gilbert bound~\eq{2-GV} to constant-weight
codes. It is shown in~\cite{Levenshtein71} 
that\vspace{-1ex}
\be{Odlyzko}
\Andw
\ \ge \
\frac{\Strut{0ex}{1.5ex}|\Fnw|}{\Strut{2.5ex}{0ex} V(n,d{-}1,w)}
\ = \
\frac{\displaystyle {n \choose w}}
     {\displaystyle \sum_{i=0}^{d-1}{w \choose i}{n{-}w \choose i}}
\ee
where $\Fnw$ is the set of binary vectors of length $n$ and weight $w$
and $V(n,d,w)$ 
is the volume of a sphere of radius~$d$\pagebreak[3.99]
in the Johnson metric.
Using the same approach as in Theorems \ref{exact} and \ref{q-ary},
we can improve upon the bound in~\eq{Odlyzko} as follows.\vspace{-1.25ex}

\begin{theorem}
\label{Andw}
Let $n$, $d$, and $w$ be three
positive integers such that $d \le w \le n/2$.
For positive~integers $i$, $j$, $k$, all less than~or~equal to $w$,
define $p^k_{i,j}$ as follows
\be{pijk}
p^k_{i,j}
\ \ \deff\ \
\sum_{l=a}^{b}
{n{-}w{-}k \choose l}
{k \choose i{-}l}
{k \choose j{-}l}
{w-k \choose i{+}j{-}k{-}l}
\ee
for all $k \le i+j$,
where 
\begin{eqnarray*}
a 
& \hspace{-.75ex}\deff\hspace{-.75ex} &
\max\{0,i-k,j-k,i+j-w\}\\
b 
& \hspace{-.75ex}\deff\hspace{-.75ex} &
\min\{i,j,i+j-k,n-w-k\}
\end{eqnarray*}
Set $p^k_{i,j} = 0$\, for $k > i+j$, and define the
following quantity
\be{endw}
e(n,d,w) 
\ \ \deff\ \
\frac{1}{6}
\sum_{i=1}^{d}
\sum_{j=1}^{d} 
\sum_{k=1}^{d}
{w \choose k} 
{n{-}w \choose k}\, p^k_{i,j}
\vspace{-1.0ex}
\ee
Then\vspace{-0.50ex}
\begin{eqnarray*}
\Andw 
& \hspace{-1.5ex}\ge\hspace{-1.5ex} &
\frac{\Strut{0ex}{1.0ex} |\Fnw|} 
{\Strut{2.25ex}{0ex} V(n,d{-}1,w)}
\\[1ex]
& & \hspace{2ex}\cdot \
\frac{\log_2\! V(n,d{-}1,w)\,-\,\log_2\!\sqrt{e(n,d{-}1,w)}}{10}
\end{eqnarray*}\vspace{-1ex}
\end{theorem}

\Proof 
The underlying ``Gilbert graph'' $\cG$ can be defined as follows:
$V(\cG) = \Fnw$ and 
$\{u,v\} \in E(\cG)$ if and only if\, $2 \le d(u,v) \le 2d'$.
Now fix a vertex $z \in V(\cG)$ and consider~the
graph $\GS$ that is induced in $\cG$ by the neighborhood $N(z)$.
Clear\-ly, all such graphs are isomorphic. The numbers $p^k_{i,j}$ 
in~\eq{pijk} are precisely the intersection numbers of the 
Johnson association scheme~\cite[p.\,665]{MWS}. It follows that if $v$
is a vertex of $\GS$ such that $d(z,v) = 2k$, then the degree 
of $v$ in $\GS$ is given by
$$
\deg(v)
\ = \
\sum_{i=1}^{d'}\!\sum_{j=1}^{d'} p^k_{i,j}
$$
Hence $e(\GS) = 3e(n,d{-}1,w)$,~where~$e(n,d,w)$ is the
quantity defined in~\eq{endw}. The desired bound on $\Andw$ now follows,
as before, from \Tref{regular}.
\qed\vspace{1ex}

\looseness=-1
We leave the asymptotic analysis of~\Tref{Andw} as
an open problem for future research. 

In the original version of this paper, we have
also suggested the following problem: generalize
the results of 
\Tref{exact} and \Tref{asymptotic} 
to lattices and sphere-packings, where the counterpart 
of the \GV\ bound is the Minkowski-Hlawka 
theorem~\cite{CS,Loeliger}.
This problem was recently solved in~\cite{KLV}.
Specifically, Krivelevich, Litsyn, and Vardy~\cite{KLV}
show that using graph-theoretic methods similar to 
those of Section\,III, the classical Minkowski
bound~\cite{Minkowski} on the density of sphere
packings in $\R^n$ can be also improved by a factor 
that is linear in $n$.

Other interesting problems for
future work would be the extension of 
Theorems~\ref{exact}~and~\ref{asymptotic} 
to spherical codes~\cite{Wyner}, 
to covering codes~\cite[Section\,12.1]{CHLL}, to codes
correcting arbitrary error pat- terns~\cite{Loeliger}, 
to runlength-limited codes~\cite{KK}, and to 
more general constrained systems~\cite{MR}.
The general approach introduced in this paper
should work~whenever an underlying 
``Gilbert graph'' can be defined, and happens
to be locally sparse.

\looseness=-1
Our results herein have applications outside of coding theory 
as well. For example, the following problem arises in the study
of scalability of optical networks~\cite{PWTD}. Let $\cH_n$
be the $n$-dimen\-sional hypercube, defined in Section\,III. What
is the minimum number $\chi_d(n)$ of colors needed to color
the vertices of $\cH_n$ so that vertices at distance $\le d$
from each other have different colors? Ngo, Du, and Graham~\cite{NDG}
have recently established the following~bound
\begin{eqnarray}
\chi_d(n)
& \hspace{-1.00ex}\le\hspace{-1.0ex} &
\hspace{.25ex}2^{\trunc{\log_2\! V(n-1,d-1)} + 1}
\nonumber\\[.25ex]
\label{NDG}
& \hspace{-1.0ex}=\hspace{-1.0ex} &
\frac{2}{\displaystyle 2^{\{\log_2\! V(n-1,d-1)\}}}\
V(n{-}1,d{-}1)
\end{eqnarray}
In fact, this follows immediately from the Varshamov bound~\eq{Varshamov-f},
since given any {\sl linear\/} binary code $\C$, assigning different 
colors to the cosets of $\C$ in $\Fn$ produces a valid coloring.
While Theorems \ref{exact} and \ref{asymptotic} improve upon~\eq{Varshamov-f},
unfortunately, we do not know whether there exist {\sl linear\/} codes that 
attain \eq{th1} and/or~\eq{th2}. Nevertheless, we can 
still improve upon~\eq{NDG}, as follows: if $d/n \le 0.499$, then 
there exists a positive constant $c$\, such that
\be{NDG-improved}
\chi_d(n)
\ \le \ 
c\ 
\frac{V(n,d)}{\log_2\! V(n,d)}
\ee
\looseness=-1
This uses a result of Alon, Krivelevich, and Sudakov~\cite{AlKS},
who show that locally sparse graphs with maximum degree $\Delta$
can be colored using $O(\Delta/\log \Delta)$ colors. Specifically,
let~$G$~be~a~graph 
with maximum degree $\Delta$ such that the 
number of edges in the sub\-graphs induced in $G$ by the neighborhood
of any vertex is at most $\Delta^2/f$. Then it is proved 
in~\cite[Theorem\,1.1]{AlKS} that the chromatic number
$\chi(G)$ of $G$ satisfies $\chi(G) \le c_1\, \Delta/\log_2\! f$
for some positive constant $c_1$.
Since the Gilbert graph $\GG$, defined in Section\,III, is
$\cH_n$ to the power $(d-1)$, it should be clear that 
$\chi_{d-1}(n) = \chi(\GG)$. The Gilbert graph $\GG$ has maximum
degree $\Delta = V(n,d{-}1) - 1$. Moreover, we have shown in the previous
section that this graph is locally sparse: if $\GS$ is the graph 
induced in $\GG$ by the neighborhood of any vertex $v \in V(\GG)$,
then $e(\GS) \le c_2\, \Delta^2 / V(n,d{-}1)^{\varepsilon}$
for $\varepsilon = 0.000001$ and some positive constant $c_2$,
provided $d/n \le 0.499$. Combining this with 
\cite[Theorem\,1.1]{AlKS} establishes~\eq{NDG-improved}.
\vspace{1ex}

Finally, we would like to point out some questions concerning
Theorems \ref{exact} and \ref{asymptotic} that remain open.
Our proof of these theorems gives no hint of linearity.
Nevertheless, we ask: are there linear codes whose parameters
satisfy~\eq{th2}? It is conceivable that a suitable modification
of the Varshamov~\cite{Varshamov} argument for constructing
a parity-check matrix could produce such codes. It is 
well known that a random (linear) code meets the \GV\
bound with probability approaching $1$ as $n \to \infty$.
Thus we ask: do random codes also meet the improved version of this
bound in \Tref{asymptotic}? 
Progress on this question was recently reported by Cohen~\cite{Cohen}.
Of course, the most interesting question of all is whether
the term $\log n$ in~\eq{GV-revised} can be further
improved to a~linear term.
In other words, it it true that the \GV\ bound 
on the {rate\/} of binary codes is asymptotically exact?

%
%
\section*{Acknowledgment}

We are grateful to Alexander Barg for helpful discussions
about~\cite{BGS}.
We would also like to thank G\'erard Cohen, Ralf Koetter, Van Vu, and 
the anonymous referees for valuable comments.\vspace{1.15ex}


\end{document}